\documentclass[a4paper,11pt]{article}

\usepackage{natbib}
\usepackage[english]{babel}
\usepackage{amsmath}
\usepackage{amssymb}
\usepackage{fancyhdr}
\usepackage[latin1]{inputenc}
\usepackage{stmaryrd}
\usepackage{amsthm}
\usepackage{geometry}
\usepackage{enumerate}
\usepackage{hyperref}
\usepackage{graphicx}

\newtheorem{thm}{Theorem}[]

\newtheorem{lm}{Lemma}[section]
\newcommand{\e}{\varepsilon}
\newcommand{\E}{\mathbb{E}}
\newcommand{\dd}{\mathrm{d}}

\newcommand{\OO}{\mathrm{O}}
\newcommand{\oo}{\mathrm{o}}

\newcommand{\tr}{\mathrm{trace}}

\newenvironment{disarray}%
{\everymath{\displaystyle\everymath{}}\array}%
{\endarray} 
\makeatletter
\renewcommand{\thefigure}{\ifnum \c@section>\z@ \thesection.\fi
 \@arabic\c@figure}
\@addtoreset{figure}{section}
\makeatother

\newcommand{\auteur}[4]{ \author{ { #1}\footnote{#2} \vspace{0.5cm} 
 \\  
 \hspace{-2.0cm} \small \noindent\begin{tabular}{l}   #3 \\ #4\end{tabular}}}

\auteur{Fanny Godet}{fanny.godet@math.univ-nantes.fr}{Laboratoire de Mathématiques Jean Leray, CNRS 6629} {Université de Nantes \\
UFR Sciences et Techniques \\
2 rue de la Houssinière - BP 92208 \\
F-44322 Nantes Cedex 3 }
\title{Prediction of long memory processes on same-realisation} 

\begin{document}
\maketitle

\begin{abstract}
For the class of stationary Gaussian long memory processes, we study
some properties of the least-squares predictor of $X_{n+1}$ based on
$(X_n, \ldots, X_1)$. The predictor is obtained by projecting $X_{n+1}$
onto the finite past and the coefficients of the predictor are estimated
on the same realisation. First we prove moment bounds for the inverse
of the empirical covariance matrix. Then we deduce an asymptotic
expression of the mean-squared error. In particular we give a relation
between the number of terms used to estimate the coefficients and the
number of past terms used for prediction, which ensures the
$L^2$-sense convergence of the predictor. Finally we prove a
central limit theorem when our predictor converges to the best linear predictor based on
all the past.

\end{abstract}


\section{Introduction}
Consider $(X_t)_{t \in \mathbb{Z}}$ a stationary process with zero mean and finite variance. We wish to predict $X_{n+1}$ from the observed past $(X_1,\ldots,X_n)$ using a linear predictor i.e. a linear combination of the observed data. First we define the coefficients of the optimal predictor in the least squares sense assuming that the covariance function is known. Then we need to estimate the replace coefficients. This second step is often realised under the following restrictive hypothesis: we predict another future independent series with exactly the same probabilistic structure; the observed series is only used to compute the forecast coefficients (see for example \cite{bhansali}, \cite{lewis} or \cite{fanny}). This assumption makes the mathematical analysis easier since the prediction problem can be reduced to an estimation problem of the forecast coefficients by conditioning on the process, which we forecast. But the practitioner rarely has two independent series: one to estimate the model, one to predict. He has to estimate the forecast coefficients on the same realisation as the forecast one. In the following we concentrate on this case called same-realisation prediction.\\
The performance of the predictor depends on two parameters: the dimension of the subspace on which we project and the number of available data to estimate the forecast coefficients. To reduce the prediction error, it is reasonable to increase the dimension of the space onto which we project, as more and more observations become available. But when the dimension and then the number of forecast coefficients increase, the estimation of these coefficients becomes more difficult and can affect the mean-squared error. \\
 When the spectral density of the process $(X_t)_{t \in \mathbb{Z}}$ exists, is bounded and bounded away from 0 (this is typical of the short memory case) \cite{ing} and \cite{kunitomo} have studied the mean-squared prediction error for same-realisation prediction. 
  The mean squared error for same-realisation prediction can be approximated by the sum of two terms: one due to the goodness of fit and one due to the model complexity. In the short-memory case, it is interesting to remark that the approximations of the mean-squared error for same and independent realisation prediction are the same. \\The performance of the least-squares predictor of long memory time series is still left unanswered. And in this case, the asymptotic equivalence between the mean-squared error in same and independent realisation should not be taken for granted since the autocovariance function decays more slowly than in the short memory case.\\ 
 The paper is organised as follows. In Sections \ref{mom} and \ref{erreur} we generalise the results of \cite{ing} to find an asymptotic expression of the mean-squared error for long memory time series. The mean-squared error is approximated by the same function as in the short memory case but under more restrictive conditions on the number of available observations and on the model complexity. In the last section, we prove a central limit theorem. More precisely, we prove the convergence in distribution of the normalised difference between our predictor and the Wiener-Kolmogorov predictor, which is the least-squares predictor knowing all the past. The normalisation is different from the short memory case since it is given by the goodness of fit of the projection.\\
\paragraph{Definition of the Predictor} 
\par Let $(X_t)_{t \in \mathbb{Z}}$ a stationary process with zero mean and finite variance. We assume that the autocovariance function $\sigma$ of the process is known. Our goal is to predict $X_{n+1}$, using the $k$ previous observed data. The optimal linear predictor is defined as the projection mapping onto the closed span of the subset $\{X_n,\ldots,X_{n-k+1}\}$ in the Hilbert space L$^2(\Omega,\mathcal{F},\mathbb{P})$ with inner product $<X,Y>=\E(X'Y)$ where $X'$ denotes the transpose of the vector $X$. It is the least-squares predictor knowing $(X_{n-k+1}, \ldots ,X_n)$. We denote by $\widetilde{X}_{n+1}(k)$ this predictor and by $-a_{j,k}$ the theoretical prediction coefficients i.e.:\begin{equation}
\widetilde{X}_{n+1}(k)= \sum_{j=1}^k\left( -a_{j,k}\right) X_{n+1-j}.\label{pred}
\end{equation}
They are given by (see \cite{coefftronccourtemem} Section 5.1):
\begin{equation}
 \left( \begin{array}{c}
 a_{1,k}\\
 \vdots \\
 a_{k,k}
\end{array}\right) 
=-\Sigma(k)^{-1} \left( \begin{array}{c}
\sigma(1)\\
 \vdots \\
\sigma(k)
\end{array}\right)\label{defajk}
\end{equation}
where  $\Sigma(k)$ is the covariance matrix of the vector $(X_1,\ldots,X_k)$.
 \paragraph*{Estimation of the Forecast Coefficients}
 When the autocovariance function $\sigma$ of the process $(X_t)_{t \in \mathbb{Z}}$ is unknown, we can plug-in an estimate of the prediction coefficients $(-a_{j,k})$ in \eqref{pred}. The estimate is constructed from the last $n$ observations $(X_n, \ldots, X_1)$ and our predictor is the projection of the last $k$ observations ($k\leq n$). The covariance matrix is estimated by:
\begin{equation}
\widehat{\Sigma}_n(k):=\frac{1}{n-K_n+1}\sum_{j=K_n}^n\textbf{X}_j(k)\textbf{X}_j'(k) \label{defmatcov}
\end{equation}
where  \begin{equation}\textbf{X}_j'(k):=(X_j,\ldots,X_{j-k+1})\label{defX}\end{equation}
and where $K_n$ is the maximum dimension of the subspace, onto which we project i.e. we will study the family of predictors $(\widetilde{X}_{n+1}(k))_{1\leq k \leq K_n}$. $K_n$ will be an increasing sequence of integers which can be bounded or can go to infinity. 
\\
The prediction coefficients $a_{i,k}$ are estimated from $(X_1, \ldots ,X_n)$ by:
$$\widehat{\textbf{a}}'(k)=(-\widehat{a}_{1,k}, \ldots , -\widehat{a}_{k,k})=\widehat{\Sigma}^{-1}_n(k)\frac{1}{n-K_n+1}\sum_{j=K_n}^{n-1}\textbf{X}_j(k)X_{j+1}.$$
The resulting one-step predictor is:
\begin{equation}
  \widehat{X}_{n+1}(k)=\textbf{X}_n'(k)\widehat{\textbf{a}}(k) \label{defpred}
\end{equation}
In this paper we use $C$ to denote generic positive constants that are independent of the sample size $n$ but may depend on the distributional properties of the process $(X_t)_{t \in \mathbb{Z}}$. Moreover $C$ may also stand for different values in different equations.\\
The following assumptions on the process $(X_t)_{t \in \mathbb{Z}}$ are essential to the results presented in the paper. There exists $d \in ]0,1/2[$ such that:
\begin{itemize}
\item[H.1] The stationary process $(X_t)_{t \in \mathbb{Z}}$ is Gaussian and admits an infinite moving average representation and an infinite autoregressive representation as follows:
\begin{equation}
\e_t =\sum_{j=0}^{+ \infty} a_j X_{t-j}\: \:\:\:\textrm{and}\:\:\:\:X_t =\sum_{j=0}^{+ \infty} b_j \e_{t-j} \label{rep}
\end{equation}
with $a_0=b_0=1$, for any $j\geq 1$ and for any $\delta>0$, $\vert a_j\vert\leq Cj^{-d-1+ \delta}$ and $\vert b_j \vert\leq Cj^{d-1+ \delta}$ and $(\e_t)_{t \in \mathbb{Z}}$ is a white noise process. These assumptions on the coefficients are verified by both long memory and short memory processes;
\item[H.2] The covariance $\sigma(k)$ is equivalent to $ L(k) k^{2d-1}$ as $k$ goes to infinity, where $L$ is a slowly varying function (i.e. for every $\alpha>0$, $x^\alpha L(x)$ is ultimately increasing and $x^{-\alpha} L(x)$ is ultimately decreasing). Under this assumption the autocovariances are not absolutely summable and thus the process is long memory process; 
\item[H.3] The spectral density of the process $(X_t)_{t \in \mathbb{Z}}$ exists and has a strictly positive lower bound;
\item[H.4] The coefficients $(a_j)_{j \in \mathbb{N}}$ verify:
 \begin{equation}
 a_j \underset{j \rightarrow + \infty}\sim L(j) j^{-d-1} \label{equivaj}
\end{equation}
with $L$ a slowly varying function.
\end{itemize}

For example, the assumptions H.1-H.4 hold for the most studied long memory process, the Gaussian FARIMA process, which is the stationary solution to the difference equations:
\begin{equation}
\phi(B) (1-B)^d X_n = \theta(B) \e_n \label{farima}
\end{equation}
where $(\e_n)_{n \in \mathbb{Z}}$ is a white noise series with mean zero, $B$ is the backward shift operator and $\phi$ and $\theta$ are polynomials with no zeroes in the unit disk.\\
We only use assumptions H.1-H.3 to give an asymptotic expression of the mean-squared error of the predictor. Assumption H.4 is a more restrictive assumption used to prove a central limit theorem for our predictor. \\
Assumption H.2 does not imply the bound on the coefficients $(a_j)_{j \in \mathbb{N}}$ and $(b_j)_{j \in \mathbb{N}}$ given in assumption H.1. \cite{inouepartialautocorr} has proved that the asymptotic expression of the autocovariance $\sigma(k) \sim L(k) k^{2d-1}$ implies:
\begin{displaymath}
 b_j \sim j^{d-1} \sqrt{\frac{L(j)}{B(d,1-2d)}}\; \textrm{ as } j\rightarrow + \infty
\end{displaymath}
and 
\begin{displaymath}
 a_j\sim j^{-d-1} \left( \sqrt{\frac{L(j)}{B(d,1-2d)}}\right) ^{-1} \frac{d \sin(\pi d)}{\pi}\; \textrm{ as } j\rightarrow + \infty
\end{displaymath}
if we assume that the sequences $(b_j)_{j \in \mathbb{N}}$ and $(a_j)_{j \in \mathbb{N}}$ are eventually decreasing to zero and $b_j \geq 0$ for all $j \in \mathbb{N}$. Such assumptions on the sign of the sequence $(b_j)_{j \in \mathbb{N}}$ or its monotonicity are not necessary for example to prove Lemma \ref{lm0} and to find moment bounds for the inverse sample covariance matrix. \\

\section{Moment bounds} \label{mom}
In this section, we establish moment bounds for the inverse sample covariance matrix  and apply these results to obtain the rate of convergence of $\widehat{\Sigma}_n(k)$ to $\Sigma(k)$. \\
 Throughout the paper, $\lambda_{min}(Y)$ and  $\lambda_{max}(Y)$ are respectively the smallest and the largest eigenvalues of the matrix $Y$. We equip the set of matrices with the norm \begin{equation}\Vert Y\Vert^2=\lambda_{max}(Y'Y) \label{norme}\end{equation}
 (see for example \cite{inegalite}). For a symmetric matrix, this norm is equal to the spectral radius and for a vector $(X_1,\ldots,X_n)$, it is equal to $\sqrt{\sum_{i=1}^n X_i^2}$. This norm is a matrix norm that verifies: for any matrices $A$ and $B$:
 \begin{equation}
 \Vert AB \Vert \leq\Vert A\Vert\Vert B\Vert\label{normematricielle}.
 \end{equation}

\begin{lm}\label{lm0}
Let $(K_n)_{n \in \mathbb{N}}$ an increasing sequence of positive integers satisfying $K_n=\oo(\sqrt{n})$. Assume (H.1).
Then, for any $q>0$, for any $\theta>0$ and for any $1 \leq k \leq K_n$,
\begin{displaymath}
\E\left[ \lambda_{min}^{-q}\left( \widehat{\Sigma}_n(k)\right) \right] = \OO\left(k^{(2+\theta)q} \right) 
\end{displaymath}
where $\widehat{\Sigma}_n(k)$ is defined in \eqref{defmatcov}.
\end{lm}
\begin{proof}The sketch of the proof is the same as that of Lemma 1 of \cite{ing}. The arguments are the following:
\begin{enumerate}
 \item the series $\sum_{j=1}^{+ \infty}|a_j|$ converges;
 \item the cumulative distribution function of the random variable $\e_t$ is a Lipschitz function and we may choose a Lipschitz constant independent of $t$. For any integer $t$ and for any reals $x$ and $y$, there exists $C$ independent of $t$ such that:
 \begin{displaymath}
 \vert\mathbb{P}( \e_t<x)-\mathbb{P}( \e_t<y)\vert \leq C \vert x-y\vert.
\end{displaymath}

\end{enumerate}
In our context these two conditions are satisfied. The sequence $(a_j)_{j \in \mathbb{N}}$ is summable under assumption H.1. Since we have assumed that the process $(X_t)_{t \in \mathbb{Z}}$ is Gaussian, $(\e_t)_{t \in \mathbb{Z}}$ is a sequence of independent and identically distributed Gaussian random variables. The distribution function of the process $\e_t$ is independent of $t$ and is a Lipschitz function. 
\end{proof}
For $n$ sufficiently large Lemma \ref{lm0} guarantees that $\widehat{\Sigma}^{-1}_n(k)$ almost surely exists as the minimum eigenvalue of $\widehat{\Sigma}_n(k)$ is almost surely positive. We also obtain an upper bound for the mean of the maximum eigenvalue of $\widehat{\Sigma}^{-1}_n(k)$. But this upper bound is not uniform as $k \rightarrow + \infty$ and therefore does not provide an asymptotic equivalent of the prediction error. Nevertheless the bound given in Lemma \ref{lm0} is a the basis of the following theorem.
\begin{thm} \label{th1}
Assume that the process $(X_t)_{t \in \mathbb{Z}}$ verifies the hypotheses H.1-H.3 
\begin{itemize}
\item if $d \in]0,1/4[$ and if there exists $\delta>0$ such that $K_n^{2+\delta}=\OO(n)$ then for all $q>0$ and for all $1 \leq k \leq K_n$:
\begin{equation}
\E \Vert \widehat{\Sigma}^{-1}_n(k)\Vert^q=\OO(1) \label{majmoment}
\end{equation}
and
\begin{equation}
 \E \Vert \widehat{\Sigma}^{-1}_n(k) - \Sigma^{-1}(k)\Vert^{q/2} \leq C\left(  \frac{K_n^2}{n-K_n+1}\right) ^{q/4} \label{majdiff}
\end{equation}
for sufficiently large $n$ ;
\item if $d \in]1/4,1/2[$ and if there exists $\delta>0$ and $\delta'>0$ such that $K_n^{2+\delta}=\OO(n^{2-4d-\delta'})$ then for all $q>0$ and for all $1 \leq k \leq K_n$:
\begin{equation}
\E \Vert \widehat{\Sigma}^{-1}_n(k)\Vert^q=\OO(1)\label{majmomentbis}
\end{equation}
and
\begin{equation}
 \E \Vert \widehat{\Sigma}^{-1}_n(k) - \Sigma^{-1}(k)\Vert^{q/2} \leq C\left(  \frac{K_n^2L^2(n-K_n+1)}{(n-K_n+1)^{2-4d}}\right) ^{q/4} \label{majdiffbis}
\end{equation}
for sufficiently large $n$;
\item    if $d =1/4$ and if there exists $\delta>0$ and $\delta'>0$ such that $K_n^{2+\delta}=\OO(n^{1-\delta'})$ then for all $q>0$ and for all $1 \leq k \leq K_n$:
\begin{equation}
\E \Vert \widehat{\Sigma}^{-1}_n(k)\Vert^q=\OO(1)\label{majmomentter}
\end{equation}
and
\begin{equation}
 \E \Vert \widehat{\Sigma}^{-1}_n(k) - \Sigma^{-1}(k)\Vert^{q/2} \leq C\left(  \frac{K_n^2L^2(n-K_n+1)\log(n-K_n+1)}{(n-K_n+1)}\right) ^{q/4} \label{majdiffter}
\end{equation}
for sufficiently large $n$.                                                       \end{itemize}
\end{thm}

In the proof of Theorem \ref{th1}, we need the following lemma.
\begin{lm} \label{lm0bis}If the process $(X_t)_{t \in \mathbb{Z}}$ verifies (H.2), if $1\leq k \leq K_n$ and \begin{itemize} 
 \item if $d\in ]0,1/4[$, then for all $q>0$,                \begin{equation}\label{approxmatcov}
\E \Vert \widehat{\Sigma}_n(k)-\Sigma(k) \Vert^q \leq C\left( \frac{K_n^2}{n-K_n+1}\right)^\frac{q}{2};
\end{equation}
\item if $d\in ]1/4,1/2[$, then for all $q>0$,
\begin{equation}\label{approxmatcov2}
\E \Vert \widehat{\Sigma}_n(k)-\Sigma(k) \Vert^q \leq C\left( \frac{K_n^2L^2(n-K_n+1)}{(n-K_n+1)^{2-4d}}\right)^\frac{q}{2};
\end{equation}
\item if $d=1/4$, then for all $q>0$,
\begin{equation}\label{approxmatcov3}
\E \Vert \widehat{\Sigma}_n(k)-\Sigma(k) \Vert^q \leq C\left( \frac{K_n^2L^2(n-K_n+1)\log(n-K_n+1)}{(n-K_n+1)}\right)^\frac{q}{2}.
\end{equation}
\end{itemize}
\end{lm}
\begin{proof}We only prove the inequalities \eqref{approxmatcov}, \eqref{approxmatcov2}  and \eqref{approxmatcov3} for $q>2$. The general case ($q>0$) easily follows from Jensen's inequality. We consider the matrix norm $\Vert .\Vert_E$ (see \cite{ciarlet}) defined for all matrix  $Y=(y_{i,j})_{1 \leq i,j \leq k}$ by
\begin{displaymath}
 \Vert Y\Vert_E =\sqrt{\sum_{i=1}^{k}\sum_{j=1}^{k} y_{i,j}^2}.
\end{displaymath}
Since the matrix $\widehat{\Sigma}_n(k)-\Sigma(k)$ is symmetric, we have
\begin{displaymath}
\Vert \widehat{\Sigma}_n(k)-\Sigma(k)\Vert \leq \Vert \widehat{\Sigma}_n(k)-\Sigma(k)\Vert_E.
\end{displaymath}
We obtain:
\begin{eqnarray}
 \Vert \widehat{\Sigma}_n(k)-\Sigma(k)\Vert^q &\leq &  \Vert\widehat{\Sigma}_n(k)-\Sigma(k)\Vert_E^q \nonumber\\
& \leq & \left(\sum_{i=1}^{k}\sum_{j=1}^{k} \left( \hat{\sigma}_{i,j}-\sigma(i-j)\right)^2 \right)^{q/2} \label {Jensen1}
\end{eqnarray}
where $\hat{\sigma}_{i,j}$ and $\sigma(i-j)$ denote respectively the $(i,j)$ entries of the matrices $\widehat{\Sigma}_n(k)$ and $\Sigma(k)$. \\
Applying Jensen's inequality to \eqref{Jensen1} because $q/2>1$, we have:
\begin{displaymath}
\Vert \widehat{\Sigma}_n(k)-\Sigma(k)\Vert^q \leq\frac{k^q}{k^2}\sum_{i=1}^{k}\sum_{j=1}^{k}\vert \hat{\sigma}_{i,j}-\sigma(i-j)\vert^q.
\end{displaymath}
It follows that:
\begin{equation}
\E \Vert \widehat{\Sigma}_n(k)-\Sigma(k)\Vert^q \leq k^{q-2}\sum_{i=1}^{k}\sum_{j=1}^{k}\E \vert \hat{\sigma}_{i,j}-\sigma(i-j)\vert^q \label{premmaj}
\end{equation}
 Now we derive the limiting distribution of $\hat{\sigma}_{i,j}-\sigma(i-j)$ to find an asymptotic expression of $\E \vert \hat{\sigma}_{i,j}~-~\sigma(i-j)\vert^q$. We shall work with the definition of the empirical covariances. By \eqref{defmatcov}, we have:
\begin{equation}
\hat{\sigma}_{i,j}= \frac{1}{n-K_n+1}\sum_{l=K_n}^n X_{l+i-1} X_{l+j-1} \underset{\mathcal{L}}=\frac{1}{n-K_n+1}\sum_{l=1}^{n-K_n+1}X_l X_{l+j-i},\label{redefcovemp}
\end{equation}
where the second equality is ensured by the strict stationarity of the process. Without loss of generality, we assume $j\geq i$. The right term of \eqref{redefcovemp} can be written:
\begin{equation}
\frac{1}{n-K_n+1}\sum_{l=1}^{n-K_n+1}X_l X_{l+j-i} =\frac{1}{n-K_n+1}\textbf{X}'_1(n-K_n+1+j-i) T_{i,j}\textbf{X}_1(n-K_n+1+j-i) \label{apptaqqu}
\end{equation}
where $\textbf{X}_1(n-K_n+1+j-i)$ is defined in \eqref{defX} and the entries of the matrix $T_{i,j}$ verify $$t_{i,j}(s,t)=
\begin{cases}
1/2 & \textrm{if}\:  \vert s-t\vert=j-i\\
0 & \textrm{otherwise}.
\end{cases}$$
$T_{i,j}$ is a Toeplitz matrix because it has symbol $g_{i,j}(x)=\cos\left( \left(j-i \right)x \right) $ i.e. $t_{i,j}(s,t)=\int_{-\pi}^{\pi}g_{i,j}(x) \cos(t-s)x)\dd x$.
\par Under Assumption H.2 with $d\in ]0,1/4[$, the spectral density verifies in a neighbourhood of~0:
\begin{displaymath}
 f(x)=\OO\left( x^{-1/2}\right)  
\end{displaymath}(see \cite{zygmund} Chap. 5 Theorem 2.6).
By applying Theorem 2 of \cite{centrallimit} to \eqref{apptaqqu}, we obtain the following convergence:  
\begin{equation}
\frac{(n-K_n+1)\left( \hat{\sigma}_{i,j}-\sigma(i-j)\right) }{\sqrt{n-K_n+1+j-i}} \underset{n\rightarrow +\infty}\Longrightarrow \mathcal{N}\left(0, 4 \pi\int_{-\pi}^{\pi} f^2(\lambda)\cos^2\left(\left( i-j\right) \lambda \right) \dd \lambda  \right) \label{cvloi}
\end{equation}
where $\Longrightarrow$ denotes the convergence in distribution.
This convergence in distribution follows from the convergence of all the cross-cumulants and hence the convergence of all the moments of the left term of \eqref{cvloi}. From this convergence in distribution we can deduce an asymptotic expression of the moments of $\hat{\sigma}_{i,j}-\sigma(i-j)$. If $q$ is even, we have an asymptotic equivalent as $n \rightarrow + \infty$:
\begin{equation}
\E \vert \hat{\sigma}_{i,j}-\sigma(i-j)\vert^q \underset{n \rightarrow +\infty} \sim \left( \frac{\sqrt{n-K_n+1+j-i}}{n-K_n+1}\right)^q\E\left[ \vert Y\vert^q\right] , \label{equivmom}
\end{equation}
where $Y$ is a Gaussian random variable which has for probability distribution the right term of \eqref{cvloi}. The $q$th-order absolute moment has the form:
\begin{equation}
E\left[ \vert Y\vert^q\right] =
\frac{q!}{2^{q/2} (q/2)!}\sigma_Y^{q}
 \label{momgauss}
\end{equation}
Moreover notice that for all $(i,j)$:
\begin{equation}
\sigma_Y \leq \sqrt{4 \pi\int_{-\pi}^{\pi} f^2(\lambda)\dd \lambda}:=M \label{defM}.
\end{equation}
Thus \eqref{equivmom}, \eqref{momgauss} and \eqref{defM} imply for sufficiently large $n$:
\begin{eqnarray*}
\E \vert \hat{\sigma}_{i,j}-\sigma(i-j)\vert^q  &\leq &
\left(  \frac{\sqrt{n-K_n+1+j-i}}{n-K_n+1}\right)^{q} \frac{q!}{(q/2)!}M^{q} \\
&\leq & C\left( \frac{1}{\sqrt{n-K_n+1}}\right)^{q}
\end{eqnarray*}
with $C$ independent of $(i,j)$. The result \eqref{approxmatcov} follows from the previous inequality and inequality \eqref{premmaj} for $q>2$. 
\par For any $d  \in]1/4,1/2[$, we apply the proposition of \cite{rosenblatt} which gives the following convergence in distribution:
\begin{displaymath}
\frac{(n-K_n+1)\left( \hat{\sigma}_{i,j}-\sigma(i-j)\right) }{L(n-K_n+1)(n-K_n+1+j-i)^{2d}} \underset{n\rightarrow +\infty}\Longrightarrow \mathcal{R}(1)
\end{displaymath}
where $\mathcal{R}$ is a Rosenblatt process. This convergence in distribution is obtained by proving the convergence of the cumulants of any order and hence the convergence of the moments of any order. This limit does not depend on the difference $(j-i)$. Similarly to the proof of \eqref{equivmom} we show for any even integer $q$:
\begin{displaymath}
\E \vert \hat{\sigma}_{i,j}-\sigma(i-j)\vert^q  \leq C\left( \frac{L(n-K_n+1)}{(n-K_n+1)^{1-2d}}\right)^{q}
\end{displaymath}
where $C$ does not depend on $(j-i)$. This inequality and \eqref{premmaj} yield the desired result.
\par Finally for $d=1/4$, we apply Theorem 4 of \cite{convcov}, which gives the following convergence in distribution:
\begin{displaymath}
 \sqrt{\frac{(n-K_n+1)}{\log(n-K_n+1)L^2(n-K_n+1)}}\left( \hat{\sigma}_{i,j}-\sigma(i-j)\right) \underset{n\rightarrow +\infty}\Longrightarrow \mathcal{N}\left(0, \sigma^2 \right) 
\end{displaymath}
where $\sigma$ depends on the process $(X_t)_{t \in \mathbb{Z}}$ but not on $(i,j)$. This convergence in distribution is obtained by proving the convergence of the cumulants of any order and hence the convergence of the moments. We then obtain that for any even integer $q$:
\begin{displaymath}
 \E \vert \hat{\sigma}_{i,j}-\sigma(i-j)\vert^q  \leq C\left(\frac{L^2(n-K_n+1)\log(n-K_n+1)}{(n-K_n+1)}\right)^{q/2}
\end{displaymath}
with $C$ independent of $(i,j)$. The result \eqref{approxmatcov3} follows from this inequality and \eqref{premmaj}. 
\end{proof}

We now prove Theorem \ref{th1}.
\begin{proof}[Proof of Theorem \ref{th1}]Since $\Vert .\Vert$ is a matrix norm (see \eqref{normematricielle}),
\begin{displaymath}
 \Vert \widehat{\Sigma}^{-1}_n(k) - \Sigma^{-1}(k)\Vert^q \leq \Vert\widehat{\Sigma}^{-1}_n(k) \Vert^q\Vert \widehat{\Sigma}_n(k) - \Sigma(k)\Vert^q\Vert \Sigma^{-1}(k)\Vert^q.
\end{displaymath}
Furthermore, by assumption H.3, the spectral density of the process $(X_t)_{t \in \mathbb{Z}}$ has a strictly positive lower bound. Thus from \cite{toeplitz}, there exists a constant $C$ such that for all $n\geq 0$:
\begin{equation}
 \Vert \Sigma^{-1}(k)\Vert^q \leq C.\label{majSn-1}
\end{equation}
Using Hölder's inequality with $1/p'+1/q'=1$, we obtain:
\begin{displaymath}
 \E\Vert \widehat{\Sigma}^{-1}_n(k) - \Sigma^{-1}(k)\Vert^q\leq C\left( \E\Vert\widehat{\Sigma}^{-1}_n(k) \Vert^{qq'}\right)^{1/q'}\left(  \E\Vert\widehat{\Sigma}_n(k) - \Sigma(k)\Vert^{qp'}\right)^{1/p'}.
\end{displaymath}
By Lemma \ref{lm0}, we have for all $\theta>0$ and  for large $n$:
\begin{displaymath}
 \left(\E\Vert\widehat{\Sigma}^{-1}_n(k) \Vert^{qq'}\right)^{1/q'} \leq C (k^{2+\theta})^q 
\end{displaymath}
then
\begin{equation}
\E\Vert \widehat{\Sigma}^{-1}_n(k) - \Sigma^{-1}(k)\Vert^q\leq C (k^{2+\theta})^q \left(  \E\Vert\widehat{\Sigma}_n(k) - \Sigma(k)\Vert^{qp'}\right)^{1/p'}\label{maj1}.
\end{equation}
We now apply Lemma \ref{lm0bis}. In order to treat together the three situations $d \in ]0,1/4[$, $d \in]1/4,1/2[$ and $d =1/4$, we define $h(n)$ by:
\begin{displaymath}
 h(n)=\left\lbrace \begin{disarray}{ll}
 \dfrac{K_n^2}{n-K_n+1}& \mathrm{if } \:d \in ]0,1/4[\\
 & \\
 \dfrac{K_n^2L^2(n-K_n+1)}{(n-K_n+1)^{2-4d}}&\mathrm{if } \:d \in]1/4,1/2[ \\
 & \\
 \dfrac{K_n^2L^2(n-K_n+1)\log(n-K_n+1)}{n-K_n+1}& \mathrm{if } \:d =1/4
\end{disarray}
\right. 
\end{displaymath}
For large $n$, we then obtain:
\begin{equation}
\left(\E\Vert\widehat{\Sigma}_n(k) - \Sigma(k)\Vert^{qp'}\right)^{1/p'} \leq C (h(n))^{q/2} \label{maj2}.
\end{equation}
From inequality \eqref{maj1} and the bound \eqref{maj2}, we obtain that there exists $ \theta>0$ such that for sufficiently large $n$:
\begin{equation}
  \E\Vert \widehat{\Sigma}^{-1}_n(k) - \Sigma^{-1}(k)\Vert^q\leq C(k^{4+\theta}h(n))^{q/2}.\label{2.24}
\end{equation}
By inequalities \eqref{2.24} and \eqref{majSn-1}, we have:
\begin{equation}
 \E\Vert\widehat{\Sigma}^{-1}_n(k) \Vert^q \leq C\left(1+ \left(k^{4+\theta}h(n)\right)^{q/2}\right) .\label{majmominv}
\end{equation}
This inequality is not sufficient to obtain \eqref{majmoment} and \eqref{majmomentbis} since under the assumptions of Theorem  \ref{th1}, $\left(k^{4+\theta}h(n)\right)^{q/2}$ is not necessarily bounded. We have to improve the intermediate inequality \eqref{majmominv}. \\
The Cauchy-Schwarz inequality and \eqref{majSn-1} give:
\begin{equation}
 \E\Vert \widehat{\Sigma}^{-1}_n(k) - \Sigma^{-1}(k)\Vert^{q/2} \leq C \left(\E\Vert\widehat{\Sigma}^{-1}_n(k) \Vert ^q \right)^{1/2} \left(\E\Vert\widehat{\Sigma}_n(k) - \Sigma(k)\Vert^q \right)^{1/2} .\label{CS}
\end{equation}
And there exists $C>0$ independent of $q$ such that:
\begin{equation}
\E\Vert\widehat{\Sigma}^{-1}_n(k) \Vert^{q/2} \leq C \left(\E\Vert\Sigma^{-1}(k)\Vert^{q/2}+\E\Vert \widehat{\Sigma}^{-1}_n(k) - \Sigma^{-1}(k)\Vert^{q/2} \right).\label{29}
\end{equation}
Inequalities \eqref{CS}, \eqref{majSn-1}, \eqref{majmominv} and Lemma \ref{lm0bis} imply that:
\begin{equation}
 \E\Vert\widehat{\Sigma}^{-1}_n(k) \Vert^{q/2} \leq C\left(1+ \left(k^{4+\theta}h(n)^2\right)^{q/4}\right).\label{30}
\end{equation}
Repeating $s-1$ times this argument (i.e. using inequalities \eqref{29}, \eqref{CS}, \eqref{majSn-1}, \eqref{30} and Lemma \ref{lm0bis}), one has for large $n$:
\begin{equation}
 \E\Vert\widehat{\Sigma}^{-1}_n(k) \Vert^{q2^{-s}}\leq C\left(1+ \left(k^{4+\theta}h(n)^{(1+s)}\right)^{q/2^{-(s+1)}}\right).\label{2.26}
\end{equation}
By assumption there exists $\delta >0$ such that $h(n)k^\delta$ converges to 0 as $n$ tends to infinity, therefore there exists $s$, such that $ \E\Vert\widehat{\Sigma}^{-1}_n(k) \Vert^{q2^{-s}}$ is bounded. Since $q$ in \eqref{2.26} is arbitrary, \eqref{majmoment} and \eqref{majmomentbis} are proved.
Inequalities \eqref{majdiff} and \eqref{majdiffbis} follow from \eqref{CS} and from Lemma \ref{lm0bis}. 
\end{proof}
In the following section, we establish an asymptotic expression for the mean-squared prediction error of the least-squares predictor using the sharp upper bound for $\E\Vert\widehat{\Sigma}^{-1}_n(k) \Vert^{q}$ given in Theorem~\ref{th1}.

\section{The mean-squared prediction error of the least-squares predictor}\label{erreur}
In this section, our goal is to give an asymptotic expression of the mean-square prediction error of the predictor defined in \eqref{defpred}. First we decompose the forecast error in:
\begin{equation}
X_{n+1}-\widehat{X} _{n+1}(k)= \e_{n+1}+f(k)+\mathcal{S}_n(k) \label{dec}
\end{equation}
where $\e_{n+1}$ is the innovation white noise at time $n+1$ and cannot be forecast, $\mathcal{S}_n(k)$ is the error due to the projection onto the closed span of the subset ${X_n, \ldots, X_{n-k+1}}$, and $f(k)$ is the error due to the estimation of the prediction coefficients. More precisely if we set $a_{i,k}=0$ for $i>k$ (for $j\leq k$, the coefficients $a_{j,k}$ are defined in \eqref{defajk}), we have
\begin{equation}
\mathcal{S}_j(k)=-\sum_{i=1}^{+\infty}(a_i-a_{i,k})X_{j+1-i} \label{defSn}
\end{equation}and with $\e_{j+1,k}$ equal to the forecast error of $X_{j+1}$ due to the projection onto $(X_j, \ldots, X_{j-k+1})$ i.e.
\begin{displaymath}
\e_{j+1,k}=X_{j+1}- P_{[X_{j-k}, \ldots ,X_j]}(X_{j+1})=X_{j+1}+\sum_{l=1}^ka_{l,k} X_{j+1-l},
\end{displaymath}
we have
\begin{displaymath}
f(k)=-\textbf{X}_n'(k)\widehat{\Sigma}^{-1}_n(k)\frac{1}{n-K_n+1}\sum_{j=K_n}^{n-1}\textbf{X}_j(k)\e_{j+1,k}
\end{displaymath}
where $\textbf{X}_n'(k)$ is defined in \eqref{defX}.

In view of \eqref{dec}, we obtain the decomposition of the mean-squared prediction error as the sum of the variance $\sigma_e^2$ of the white noise and the error due to the prediction method $\E\left( f(k)+\mathcal{S}_n(k)\right) ^2$:
\begin{displaymath}
\E\left(  X_{n+1}-\widehat{X} _{n+1}(k)\right)^2 = \sigma_e^2+\E\left( f(k)+\mathcal{S}_n(k)\right) ^2.
\end{displaymath}

\begin{thm}\label{th2}Under assumptions H.1-H.3, if we choose the sequence $(K_n)_{n \in \mathbb{N}}$ such that for some $\delta>0$:
\begin{equation}
K_n^{4}=\oo(n^{1-2d-\delta}) \label{condth2},
\end{equation}
then
\begin{displaymath}
\lim_{n \rightarrow + \infty} \max_{1 \leq k \leq K_n} \left \vert\frac{\E\left(X_{n+1} - \widehat{X}_{n+1}(k) \right)^2 -\sigma_\e^2}{L_n(k)}-1 \right \vert=0
\end{displaymath}
where
\begin{equation}
L_n(k)=\E\left( \mathcal{S}_n(k)\right)^2+\frac{k}{n-K_n+1}\sigma_\e^2 \label{defL}
\end{equation}
$\mathcal{S}_n(k)$ being defined in \eqref{defSn}.
\end{thm}
\paragraph*{Remark} If we fit a misspecified AR($k$) model to the long-memory time series $(X_n)_{n \in \mathbb{Z}}$ to forecast it, we find the same predictor as \eqref{pred}. Consequently $L_n(k)$ can be viewed as the quality of prediction by an AR model. From \eqref{defL}, this quality is the sum of the model complexity $\frac{k}{n-K_n+1}\sigma_\e^2$ and the goodness of fit $\E\left( \mathcal{S}_n(k)\right)^2$.
\begin{proof}
By \eqref{dec}, we have:
\begin{displaymath}
\left \vert\frac{\E\left( X_{n+1} - \widehat{X}_{n+1}(k) \right)^2 -\sigma_\e^2}{L_n(k)}-1 \right \vert= \left \vert\frac{E\left( f(k)+\mathcal{S}_n(k)\right)^2}{L_n(k)}-1 \right \vert.
\end{displaymath}
Our proof is divided into three steps:
\begin{enumerate}
\item we provide an approximation of $\E( f(k))^2$ which is easier to estimate. This approximation denoted by $\E( f_1(k))^2$ will be defined in \eqref{f1};
\item we show that the asymptotic equivalent of $\E( f_1(k))^2$ is $\frac{k}{n-K_n+1}\sigma_e^2$;
\item we prove that the cross-product term $\E( f(k)\mathcal{S}_n(k))$ is negligible with respect to $L_n(k)$. 
\end{enumerate}
\paragraph*{First step}
We introduce
\begin{equation}
f_1(k):=-\textbf{X}_n^{*'} (k)\Sigma^{-1}(k)\frac{1}{n-K_n+1}\sum_{j=K_n}^{n-\sqrt{n}-1}\textbf{X}_j(k)\e_{j+1},\label{f1}
\end{equation}
with
\begin{displaymath}
 \textbf{X}_n^{*'} (k)=\left(\sum_{j=0}^{\sqrt{n}/2-K_n} b_j\e_{n-j}, \ldots,\sum_{j=0}^{\sqrt{n}/2-K_n} b_j\e_{n-k+1-j}\right) .
\end{displaymath}
\begin{lm}\label{lmcalc}If the assumptions of Theorem \ref{th2} hold, then
\begin{equation}
\lim_{n \rightarrow + \infty} \max_{1 \leq k \leq K_n}\E\left(\sqrt{\frac{1}{L_n(k)}}(f(k)-f_1(k)) \right)^2 =0.\label{f-f1}
\end{equation}
\end{lm}\begin{proof}
See the appendix.
\end{proof}

\paragraph*{Second step}
We prove that
\begin{equation}
 \lim_{n \rightarrow + \infty} \max_{1 \leq k \leq K_n} \left \vert\E\left( \frac{n-K_n+1}{k\sigma_\e^2}f_1^2(k)\right) -1 \right\vert=0. \label{convf1}
\end{equation}
First observe that
\begin{eqnarray*}
 &&\E\left( \frac{n-K_n+1}{k\sigma_\e^2}f_1^2(k)\right)\\
 &=&\frac{n-K_n+1}{k\sigma_\e^2}\E\left(\textbf{X}_n^{*'} (k)\Sigma^{-1}(k)\frac{1}{n-K_n+1}\sum_{j=K_n}^{n-\sqrt{n}-1}\textbf{X}_j(k)\e_{j+1}\right)^2\\
 &=&\frac{n-K_n+1}{k\sigma_\e^2}  \E\left[  \tr \left(\textbf{X}_n^{*'} (k)\Sigma^{-1}(k)\frac{1}{n-K_n+1}\sum_{j=K_n}^{n-\sqrt{n}-1}\textbf{X}_j(k)\e_{j+1}\right)^2\right] \\
 &=&\frac{n-K_n+1}{k\sigma_\e^2}\tr\left[ \E\left(\Sigma^{-1}(k)\frac{1}{(n-K_n+1)^2}\sum_{j=K_n}^{n-\sqrt{n}-1}\textbf{X}_j(k)\e_{j+1} \sum_{l=K_n}^{n-\sqrt{n}-1}\textbf{X'}_l(k)\e_{l+1}\Sigma^{-1}(k)\textbf{X}_n^{*}\textbf{X}_n^{*'}\right) \right] 
 \end{eqnarray*}
Since the vector $\textbf{X}_n^{*'}$ and $\displaystyle{\sum_{j=K_n}^{n-\sqrt{n}-1}\textbf{X}_j(k)\e_{j+1}}$ are uncorrelated because $k \leq K_n$:
 \begin{eqnarray*}
 \E\left( \frac{n-K_n+1}{k\sigma_\e^2}f_1^2(k)\right)&=& \frac{n-K_n+1}{k\sigma_\e^2(n-K_n+1)^2}\tr \left(\Sigma^{-1}(k) (n-K_n+1-\sqrt{n})\sigma_\e^2\Sigma(k)\Sigma^{-1}(k)\Sigma^*(k)\right) \\
 &=& \tr\left(\Sigma^{-1}(k) \Sigma^*(k)k^{-1}\right) (n-K_n+1-\sqrt{n})(n-K_n+1)^{-1},
\end{eqnarray*}
where $\Sigma^*(k)$ is the covariance matrix of the vector $\textbf{X}_n^*(k)$. We note that:
\begin{displaymath}
 (n-K_n+1-\sqrt{n})(n-K_n+1)^{-1} \rightarrow 1 \;\textrm{as} \; n \rightarrow + \infty.
\end{displaymath}
So we only have to study the trace of $ \left(\Sigma^{-1}(k) \Sigma^*(k)k^{-1}\right)$. We will use the following inequality: for all $k\times k$ matrices $A$ and $B$
\begin{eqnarray*}
 \vert \tr (AB)\vert & \leq & \sqrt{\tr (AA')}\sqrt{\tr (BB')}\\
 &\leq & k \Vert A \Vert \Vert B \Vert.
\end{eqnarray*}
We obtain:
\begin{eqnarray*}
 \max_{1 \leq k \leq K_n} \left \vert \tr\left(\Sigma^{-1}(k) \Sigma^*(k)k^{-1}\right)-1 \right \vert&=& \max_{1 \leq k \leq K_n} \left \vert\tr\left(\Sigma^{-1}(k) (\Sigma^*(k)-\Sigma(k))k^{-1}\right)\right \vert\\
 &\leq &\max_{1 \leq k \leq K_n} \left \Vert\Sigma^{-1}(k) \right \Vert\left \Vert(\Sigma^*(k)-\Sigma(k))\right \Vert \\
 &\leq &\max_{1 \leq k \leq K_n} \left \Vert\Sigma^{-1}(k) \right \Vert\max_{1 \leq k \leq K_n}\left \Vert(\Sigma^*(k)-\Sigma(k))\right \Vert \end{eqnarray*} 
 $\Sigma(k)-\Sigma^*(k)$ is symmetric because $\Sigma(k)$ and $\Sigma^*(k)$ are two symmetric matrices, and its spectral norm is lower than every other matrix norm. We use the subordinate norm defined for all matrix $Y=(y_{i,j})_{1 \leq i,j\leq k}$ by:
 \begin{displaymath}
 \Vert Y\Vert_1 =\max_j \sum_{i=1}^k |y_{i,j}|.
 \end{displaymath}
 For large $n$, we obtain
 \begin{eqnarray*}
 \max_{1 \leq k \leq K_n}\Vert\Sigma^{-1}(k)\Vert\max_{1 \leq k \leq K_n}\Vert \Sigma^*(k)-\Sigma(k)\Vert&\leq &\max_{1 \leq k \leq K_n}\Vert\Sigma^{-1}(k)\Vert\max_{1 \leq k \leq K_n}\Vert\Sigma(k)-\Sigma^*(k)\Vert_1\\
 &\leq &\max_{1 \leq k \leq K_n}\Vert\Sigma^{-1}(k)\Vert\max_{1 \leq k \leq K_n}k \max_{0 \leq j \leq k-1}\sum_{l=\sqrt{n}/2-K_n+1}^{+\infty}\vert b_l b_{l+j}\vert \\
&=&\OO\left(\frac{K_n}{\left( \sqrt{n} \right)^{1-2d-\delta}} \right) 
\end{eqnarray*}
for all $\delta>0.$ Then
\begin{eqnarray*}
\max_{1 \leq k \leq K_n}\Vert\Sigma^{-1}(k)\Vert\max_{1 \leq k \leq K_n}\Vert \Sigma^*(k)-\Sigma(k)\Vert&=&\oo(1)
\end{eqnarray*}
follows from condition \eqref{condth2}.

\paragraph*{Third step}
We consider the cross-product term $\E \left( f(k) \mathcal{S}_n(k) L_n^{-1}(k)\right)$ and show that it is negligible. \cite{ing} proved that:
\begin{displaymath}
 \left |\E \left( f(k) \mathcal{S}_n(k) L_n^{-1}(k)\right) \right|= \left |\E \left( (f(k)-f_1(k)) \mathcal{S}_n(k) L_n^{-1}(k)\right)\right| . 
\end{displaymath}
By the Cauchy-Schwarz inequality:
\begin{eqnarray*}
&& \max_{1 \leq k \leq K_n}\left \vert \E \left( (f(k)-f_1(k)) \mathcal{S}_n(k) L_n^{-1}(k)\right) \right \vert \\&\leq& \left[\max_{1 \leq k \leq K_n}\E\left((f(k)-f_1(k))^2  L_n^{-1}(k)\right) \max_{1 \leq k \leq K_n}\E\left( \mathcal{S}^2_n(k)L_n^{-1}(k)\right)\right]^{1/2} .
\end{eqnarray*}
By \eqref{f-f1}, we obtain:
\begin{displaymath}
 \max_{1 \leq k \leq K_n}\E\left((f(k)-f_1(k))^2  L_n^{-1}(k)\right)=\oo(1)
\end{displaymath}
and by using the definition \eqref{defL} of $L_n(k)$, we have
\begin{displaymath}
 \max_{1 \leq k \leq K_n}\E\left( \mathcal{S}^2_n(k)L_n^{-1}(k)\right)= \OO(1).
\end{displaymath}
Finally  we have
\begin{displaymath}
 \max_{1 \leq k \leq K_n}\left \vert \E \left( (f(k)-f_1(k)) \mathcal{S}_n(k) L_n^{-1}(k)\right) \right \vert =\oo(1).
\end{displaymath}

\end{proof}
In this theorem, we have obtained an asymptotic expression of the mean squared prediction error of $\widehat{X} _{n+1}(k)$, which holds uniformly for all $1 \leq k \leq K_n$. In the short memory case i.e. assuming that the process $(X_t)_{t \in \mathbb{Z}}$ is Gaussian, admits infinite moving average and autoregressive representations defined in \eqref{rep}, that the coefficients $(a_j)_{j \in \mathbb{N}}$ verify $\sum_{j=1}^{+\infty}\sqrt{j}|a_j|< \infty$ and that the coefficients $b_j$ are absolutely summable, \cite{ing} proved that if $K_n^{2+\delta}=\OO(n)$ for some $\delta>0$:
\begin{displaymath}
\lim_{n \rightarrow + \infty} \max_{1 \leq k \leq K_n} \left \vert\frac{\E\left(X_{n+1} - \widehat{X}_{n+1}(k) \right)^2 -\sigma_\e^2}{L_n(k)}-1 \right \vert=0
\end{displaymath}
with $L_n(k)$ defined as in \eqref{defL}.\\
The term $L_n(k)$ has the same expression in the short memory case as in the long memory case. It is the sum of two terms: the first term $(k/n)\sigma_\e^2$ is proportional to the order of the model and is a measure of the complexity of the predictor, the second term $\mathcal{S}^2_n(k)$ corresponds to the goodness of fit of the model. This second term has not the same asymptotic behaviour in theshort and long memory case: for short memory time series it decays exponentially fast as a function of $k$, whereasfor long memory time series it has a Riemannian decay.\\

In the following section, we will use the proof of Theorem \ref{th2} to obtain a central limit theorem for our predictor.

\section{Central limit theorem}

Like \cite{bhansali} and \cite{lewis} for short memory processes, we search a normalisation factor to obtain a convergence in distribution of the difference between our predictor $\widehat{X} _{n+1}(K_n)$ and the Wiener-Kolmogorov predictor $\displaystyle{\widetilde{X}_{n+1}=-\sum_{j=1}^{+\infty} a_jX_{n+1-j}}$, which is the linear least-squares predictor based on all the past.

\begin{thm} \label{th3}Under assumptions H.1-H.4, if we choose the sequence $(K_n)_{n \in \mathbb{Z}}$ such that:
\begin{equation}
 K_n^4=\OO(n) \; \textrm{ and }\; K_n^{1+2d}=\oo\left( n^{1-2d}\right), \label{condth3}
\end{equation}
then 
\begin{displaymath}
  \frac{1}{\sqrt{\E[\mathcal{S}_n^2(K_n)]}}\left(\widetilde{X}_{n+1}-\widehat{X} _{n+1}(K_n) \right)\xrightarrow[n \rightarrow + \infty]{} \mathcal{N}(0,1).
\end{displaymath}

\end{thm}
\begin{proof}
The difference between our predictor $\widehat{X} _{n+1}(K_n)$ and the Wiener-Kolmogorov predictor $\widetilde{X}_{n+1}$ is equal to:
\begin{equation}
 \widetilde{X}_{n+1}-\widehat{X} _{n+1}(K_n)=f(K_n)+\mathcal{S}_n(K_n). \label{dec2}
\end{equation}
Since $(X_t)_{t \in \mathbb{Z}}$ is a Gaussian process with mean $0$, $(\sum_{i=1}^l(a_i-a_{i,K_n}) X_{t+1-i})_{t \in \mathbb{Z}}$ is a Gaussian random variable with mean $0$ for any integer $l$ . But $(\sum_{i=1}^l(a_i-a_{i,K_n}) X_{t+1-i})_{t \in \mathbb{Z}}$ converges in mean-squared sense and thus in distribution to $\mathcal{S}_n(K_n)$ as $l$ tends to infinity. Then $\mathcal{S}_n(K_n)$ is Gaussian random variable with mean 0.\\
Consequently it is enough to prove that
 \begin{displaymath}
   \frac{1}{\sqrt{\E[\mathcal{S}_n^2(K_n)]}}f(K_n) \xrightarrow[n \rightarrow + \infty]{\mathbb{P}} 0 .
 \end{displaymath}

First we search for a bound for $1/\E[\mathcal{S}_n^2(K_n)]$. 
For all integer $l$,
\begin{displaymath}
\E\left(\sum_{i=1}^l (a_i-a_{i,K_n})X_{n+1-i} \right) ^2\geq 2 \pi \underline{f} \sum_{i=1}^l(a_i-a_{i,K_n})^2
\end{displaymath}
because the spectral density $f$, which is the Toeplitz symbol of the covariance matrix, is bounded below by a positive constant $\underline{f}$ (see \cite{toeplitz}).
By taking the limit as $l \rightarrow + \infty$, we obtain:
\begin{eqnarray*}
E\left(\sum_{i=1}^{ + \infty} (a_i-a_{i,K_n})X_{j+1-i} \right) ^2&\geq&2 \pi\underline{f}  \sum_{i=1}^{+ \infty}(a_i-a_{i,K_n})^2 \\
 &\geq & 2 \pi \underline{f}\sum_{i=K_n+1}^{+ \infty}a_i^2
\end{eqnarray*}
since $a_{i,K_n}=0$ when $i>K_n$.
Under assumption H.4,  
\begin{displaymath}
\sum_{i=K_n+1}^{+ \infty}a_i^2 \underset{n \rightarrow + \infty}\sim\frac{1}{1+2d}K_n^{-2d-1}L^2(K_n)
\end{displaymath}(see Proposition 1.5.10 of \cite{bingham}).
Then for any $\delta >0$, there exists $C>0$ such that:
\begin{equation}
 \frac{1}{\E[\mathcal{S}_n^2(K_n)]}\leq C K_n^{2d+1+\delta} \label{majinvSn}
\end{equation}

 By introducing $f_1$ defined in the proof of Theorem \ref{th2}, we decompose the proof of the mean-squared convergence in two parts. We will first show that:
 \begin{equation}
 \frac{1}{\sqrt{\E[\mathcal{S}_n^2(K_n)]}}\left( f(K_n)-f_1(K_n)\right) \xrightarrow[n \rightarrow + \infty]{\mathrm{L}^2} 0 \label{annoncebis}
\end{equation} then
\begin{equation}
 \frac{1}{\sqrt{\E[\mathcal{S}_n^2(K_n)]}}f_1(K_n) \xrightarrow[n \rightarrow + \infty]{\mathrm{L}^2}0 . \label{annonce}
\end{equation}

More precisely we will prove the mean-squared convergence \eqref{annoncebis}, using the decomposition in four terms \eqref{decomp1}, \eqref{decomp2}, \eqref{decomp3} and \eqref{decomp4} of proof of Lemma \ref{lmcalc} (see appendix). Using \eqref{term1} and \eqref{majinvSn}, the term \eqref{decomp1} verifies for any $\delta>0$:
\begin{eqnarray}
  &&\E\left( \frac{1}{\sqrt{\E[\mathcal{S}_n^2(K_n)]}}\textbf{X}_n^{*'} (K_n)\Sigma^{-1}(K_n)\frac{1}{n-K_n+1}\left[ \sum_{j=K_n}^{n-\sqrt{n}-1}\textbf{X}_j(K_n)\e_{j+1}-\sum_{j=K_n}^{n-1}\textbf{X}_j(K_n)\e_{j+1}\right] \right)^2 \nonumber\\&=& \OO\left(\frac{K_n^{3+2d+\delta}}{n^{5/4}} \right). \label{anommer}
\end{eqnarray}
Under assumption \eqref{condth3}, the mean \eqref{anommer} converges to 0.\\
 Similarly for the term \eqref{decomp2} using \eqref{term2} and \eqref{majinvSn} we obtain for any $\delta>0$:
 \begin{displaymath}
 \E\left( \frac{1}{\sqrt{\E[\mathcal{S}_n^2(K_n)]}}\textbf{X}_n^{*'} (K_n)\left[\Sigma^{-1}(K_n)- \widehat{\Sigma}^{-1}_n(K_n)\right] \sum_{j=K_n}^{n-1}\textbf{X}_j(K_n)\e_{j+1}\right)^2 =\OO\left(\frac{K_n^{5+2d+\delta}}{(n-K_n+1)^2} \right) 
\end{displaymath}
which converges to 0 under assumption \eqref{condth3} for sufficiently small $\delta$. \\
For the third term \eqref{decomp3}, by \eqref{term3} and \eqref{majinvSn} we obtain:
\begin{displaymath}
 \E\left( \frac{1}{\sqrt{\E[\mathcal{S}_n^2(K_n)]}}\left[ \textbf{X}_n^{*'} (K_n)-\textbf{X}_n'(K_n)\right] \widehat{\Sigma}^{-1}_n(K_n) \sum_{j=K_n}^{n-1}\textbf{X}_j(K_n)\e_{j+1}\right)^2= \OO\left(\frac{K_n^{3+2d}}{(n-K_n+1)^{\frac{3-2d}{2}}} \right) 
\end{displaymath}
which converges to 0 under assumption \eqref{condth3}. \\
Finally the estimation of the fourth term \eqref{decomp4} is directly given in \eqref{term4}:
\begin{eqnarray*}
  &&\E\left( \sqrt{\frac{1}{\E\left( \mathcal{S}_n(K_n)\right)^2}}\textbf{X}_n' (K_n)\widehat{\Sigma}^{-1}_n(K_n)\frac{1}{n-K_n+1}\sum_{j=K_n}^{n-1}\textbf{X}_j(K_n)\left[ \e_{j+1,K_n}-\e_{j+1}\right] \right)^2 \\&=&\OO\left( \frac{K_n^{1+2d+\delta}}{(n-K_n+1)^{1-2d-\delta}}\right) 
\end{eqnarray*}
which converges to 0 under condition \eqref{condth3}.\\
Now we will prove the mean-squared convergence \eqref{annonce}.

By \eqref{convf1}:
\begin{displaymath}
\E\left( \frac{n-K_n+1}{K_n\sigma_\e^2}f_1^2(K_n)\right) \underset{n \rightarrow + \infty}\sim \frac{K_n\sigma_\e^2}{n-K_n+1}.
\end{displaymath}
Under condition \eqref{condth3}, bound \eqref{majinvSn} implies:
\begin{displaymath}
 \frac{1}{\E[\mathcal{S}_n^2(K_n)]}\frac{K_n\sigma_\e^2}{n-K_n+1} \xrightarrow[n \rightarrow + \infty]{}0.
\end{displaymath}
Then we have:
\begin{displaymath}
 \lim_{n \rightarrow + \infty} \frac{1}{\E[\mathcal{S}_n^2(K_n)]} \E\left( f_1^2(K_n)\right) =0.
\end{displaymath}

\end{proof}

\paragraph*{Remark 1}The normalisation in Theorem \ref{th3} is not an explicit function of $K_n$. Nevertheless we have a good idea of the rate of decay of $K_n^{-1}$ to $0$. We have shown in \eqref{majinvSn} that under assumption H.4 for all $\delta>0$:
\begin{displaymath}
\exists C, \: C K_n^{-2d-1-\delta} \leq \E[\mathcal{S}_n^2(K_n)] .
\end{displaymath}
In \cite{fanny2}[Theorem 3.3.1], an upper bound for the rate of convergence is proved assuming H.1-H.2: 
\begin{displaymath}
\exists C, \:\E[\mathcal{S}_n^2(K_n)] \leq C K_n^{-1}.
\end{displaymath}
For some processes, we even have an equivalent of $\E[\mathcal{S}_n^2(K_n)]$. Consider a fractionally integrated noise $(X_t)_{t \in \mathbb{Z}}$, which is the stationary solution of the difference equation:
\begin{displaymath}
 (I-B)^d X_t = \e_t
\end{displaymath}
where $(\e_t)_{t \in \mathbb{Z}}$ is a white noise with mean $0$ and constant finite variance $\sigma_\e^2$ and $B$ is the backward-shift operator. In this case, the rate of convergence is given by:
\begin{displaymath}
\exists C, \:\E[\mathcal{S}_n^2(K_n)] \sim C K_n^{-1}.
\end{displaymath}
\paragraph*{Remark 2}  In both the short and the long-memory case, the prediction error between our predictor and the Wiener-Kolmogorov predictor has the same expression $L_n(K_n)=\frac{K_n}{n}\sigma_\e^2+\E[\mathcal{S}_n^2(K_n)]$. But we do not use the same normalisation for central limit theorems.\\
In the central limit theorem for short memory processes, we only know results for independent realisation prediction i.e. when the aim is to predict an independent series which has exactly the same probabilistic structure as the observed one. \cite{bhansali} and \cite{lewis} proved a convergence in distribution of $\left(\widetilde{X}_{n+1}-\widehat{X} _{n+1}(K_n) \right) $ normalised by $\sqrt{\frac{n}{K_n\sigma^2_\e}}$ respectively in the univariate case and in the multivariate case. This normalisation corresponds to the complexity of the estimation of the projection coefficients.\\
In the long memory case the normalisation $\sqrt{\E[\mathcal{S}_n^2(K_n)] }$ is given by the rate of convergence of the predictor knowing a finite past to the linear least-squares predictor knowing the infinite past. In the long memory case the rate of convergence due to the projection decays hyperbolically and is the main term of the global error of prediction $L_n(K_n)$. On the contrary in the short memory case, the rate of convergence due to the projection decays exponentially fast and is negligible with respect to the rate of convergence due to the estimation of the projection coefficients. 

\section{Appendix}
\subsection{Preliminary lemmas}
In the following lemmas we prove subsidiary asymptotic results, which we need in the proof of Theorem \ref{th2}.
\begin{lm}\label{lem1}Assume H.2. If $q\geq 1$, then for all $\delta>0$, there exists $C$ constant such that for all $1 \leq k \leq K_n$:
\begin{equation}
 \E\left\Vert\frac{1}{\sqrt{n-K_n+1}}\sum_{j=K_n}^{n-1}\textbf{X}_j(k)(\e_{j+1,k}-\e_{j+1})\right\Vert^q \leq C\left( k (n-K_n+1)^{2d+\delta}\E\left( \mathcal{S}_n(k)\right)^2 \right) ^{q/2} \label{eqlem1}
\end{equation}
where the norm $\Vert .\Vert$ is defined in \eqref{norme}.
\end{lm}
\begin{proof}
We have
\begin{displaymath}
\frac{1}{\sqrt{n-K_n+1}}\sum_{j=K_n}^{n-1}\textbf{X}_j(k)(\e_{j+1,k}-\e_{j+1})=\frac{1}{\sqrt{n-K_n+1}}\sum_{j=K_n}^{n-1}\textbf{X}_j(k)\sum_{i=1}^{+\infty}(a_i-a_{i,k})X_{j+1-i}.
\end{displaymath}
Without loss of generality, we assume that $q>2$ since the result for $q>1$ can be obtained from the result for $q>2$ and Jensen's inequality. Observe that:
\begin{eqnarray*}
&& \left \Vert\frac{1}{\sqrt{n-K_n+1}}\sum_{j=K_n}^{n-1}\textbf{X}_j(k)\sum_{i=1}^{+\infty}(a_i-a_{i,k})X_{j+1-i}\right \Vert ^q\\
&=&(n-K_n+1)^{-q/2}
 \left( \sum_{l=0}^{k-1}\left(\sum_{j=K_n}^{n-1}X_{j-l} \sum_{i=1}^{+\infty}(a_i-a_{i,k})X_{j+1-i}  \right)^2 \right)^{q/2} 
\end{eqnarray*}
Since the function $x \mapsto x^{q/2}$ is convex on $\mathbb{R}^+$ if $q>2$, we obtain by Jensen's inequality:
\begin{displaymath}
 \left( \sum_{l=0}^{k-1}\left(\sum_{j=K_n}^{n-1}X_{j-l} \sum_{i=1}^{+\infty}(a_i-a_{i,k})X_{j+1-i}  \right)^2 \right)^{q/2} \leq  k^{-1}\sum_{l=0}^{k-1}k^{q/2}\left\vert\sum_{j=K_n}^{n-1}X_{j-l} \sum_{i=1}^{+\infty}(a_i-a_{i,k})X_{j+1-i}  \right\vert^q.
\end{displaymath}
Consequently
\begin{eqnarray}
&&\E\left\Vert\frac{1}{\sqrt{n-K_n+1}}\sum_{j=K_n}^{n-1}\textbf{X}_j(k)(\e_{j+1,k}-\e_{j+1})\right\Vert^q \nonumber \\
&\leq& k^{q/2-1}  \sum_{l=0}^{k-1}\E\left((n-K_n+1)^{-q/2} \left \vert\sum_{j=K_n}^{n-1}X_{j-l} \sum_{i=1}^{+\infty}(a_i-a_{i,k})X_{j+1-i} \right \vert^q \right). \label{similaire}
\end{eqnarray}
Furthermore
\begin{eqnarray}
\E\left(X_{j-l} \sum_{i=1}^{+\infty}(a_i-a_{i,k})X_{j+1-i} \right) &=&
\E\left(X_{j-l}\left[ \e_{j+1} -X_{j+1}-\sum_{i=1}^k a_{i,k} X_{j+1-i}\right] \right) \label{esp0}\\
&=& -\sigma_{l+1}-\sum_{i=1}^k a_{i,k} \sigma_{l+1-i}\nonumber\\
&=&0 \nonumber
\end{eqnarray}
for any integer $l \in \left[ 1 ,k \right] $ by definition of $(a_{i,k})_{1\leq i\leq k}$. Since the mean defined in \eqref{esp0} is equal to 0,
 \begin{eqnarray*}
 &&\E\left((n-K_n+1)^{-q/2} \left \vert\sum_{j=K_n}^{n-1}X_{j-l} \sum_{i=1}^{+\infty}(a_i-a_{i,k})X_{j+1-i} \right \vert^q \right)\\ &=&\E\left((n-K_n+1)^{-q/2} \left \vert\sum_{j=K_n}^{n-1}X_{j-l} \sum_{i=1}^{+\infty}(a_i-a_{i,k})X_{j+1-i}- \E\left(X_{j-l} \sum_{i=1}^{+\infty}(a_i-a_{i,k})X_{j+1-i} \right)\right \vert^q \right).
\end{eqnarray*}
And then by applying Theorem 1 of \cite{ing} to the random variable \\ $Q=\sum_{j=K_n}^{n-1}X_{j-l} \sum_{i=1}^{+\infty}(a_i-a_{i,k})X_{j+1-i}$, we obtain:
\begin{eqnarray*}
&&\E\left((n-K_n+1)^{-q/2} \left \vert\sum_{j=K_n}^{n-1}X_{j-l} \sum_{i=1}^{+\infty}(a_i-a_{i,k})X_{j+1-i} \right \vert^q \right) \\&\leq &C \left( \frac{1}{n-K_n+1}\sum_{s=K_n}^{n-1} \sum_{t=K_n}^{n-1} \sigma(s-t) \sigma^*(s-t)\right)^{q/2}, 
\end{eqnarray*}
where $\sigma^*(.)$ is the autocovariance function of the process $\left( \sum_{i=1}^{+\infty}(a_i-a_{i,k})X_{t+1-i}\right)_{t \in \mathbb{Z}} $ i.e.
\begin{displaymath}
\sigma^*(s-t)=\E\left[\left( \sum_{i=1}^{+\infty}(a_i-a_{i,k})X_{s+1-i}\right)  \left( \sum_{i=1}^{+\infty}(a_i-a_{i,k})X_{t+1-i}\right) \right] .
\end{displaymath}
As $\left \vert \sigma^*(s-t)\right \vert \leq \sigma^*(0)$,
\begin{eqnarray*}
&&\E\left((n-K_n+1)^{-q/2} \left \vert\sum_{j=K_n}^{n-1}X_{j-l} \sum_{i=1}^{+\infty}(a_i-a_{i,k})X_{j+1-i} \right \vert^q \right)\\
 & \leq &C \left( \frac{1}{n-K_n+1}\sigma^*(0)\sum_{s=K_n}^{n-1} \sum_{t=K_n}^{n-1} \vert\sigma(s-t)\vert\right)^{q/2}\\
 & \leq &C\left( \frac{1}{n-K_n+1}\sigma^*(0)\sum_{s=1}^{n-K_n+1} \sum_{t=1}^{n-K_n+1} \vert\sigma(s-t)\vert\right)^{q/2}.
\end{eqnarray*}
Under assumption H.2, we have for any $\delta >0$
\begin{displaymath}
\sum_{s=1}^{n-K_n+1} \sum_{t=1}^{n-K_n+1} \vert\sigma(s-t)\vert \leq C(n-K_n+1)^{2d+ \delta+1}.
\end{displaymath}
Then we have shown that for all $\delta >0$,
\begin{displaymath}
\E\left((n-K_n+1)^{-q/2} \left \vert\sum_{j=K_n}^{n-1}X_{j-l} \sum_{i=1}^{+\infty}(a_i-a_{i,k})X_{j+1-i} \right \vert^q \right) \leq C\left((n-K_n+1)^{2d+\delta}\sigma^*(0)\right)^{q/2}.
\end{displaymath}
Notice that:
\begin{displaymath}
\sigma^*(0)= \E\left( \mathcal{S}_n(k)\right)^2.
\end{displaymath}
And this remark allows us to conclude.
\end{proof}

\begin{lm}\label{lem2}Assume that the assumptions of Theorem \ref{th2} hold. If $q>1$, then for any $1 \leq k\leq K_n$
\begin{displaymath}
\E\left\Vert\frac{1}{\sqrt{n-K_n+1}}\sum_{j=K_n}^{n-1}\textbf{X}_j(k)\e_{j+1}\right\Vert^q \leq C k ^{q/2}
\end{displaymath}
with $C$ independent of $n$ and then of $k$.
\end{lm}
\begin{proof}
The arguments are similar to those used for verifying Lemma \ref{lem1}. Without loss of generality we assume that $q>2$, since this result and Jensen's inequality allow to conclude for $q>1$. Reasoning as for \eqref{similaire}, we have by convexity:
\begin{displaymath}
\E\left\Vert\frac{1}{\sqrt{n-K_n+1}}\sum_{j=K_n}^{n-1}\textbf{X}_j(k)\e_{j+1}\right\Vert^q \leq k^{q/2} k^{-1}\sum_{l=0}^{k-1}\E\left[ (n-K_n+1)^{-q/2} 
\left \vert 
\sum_{j=K_n}^{n-1} X_{j-l} \e_{j+1}
 \right \vert^q 
\right] .
\end{displaymath}
Applying again Theorem 1 of \cite{ing}:
\begin{displaymath}
\E\left[ (n-K_n+1)^{-q/2} \left \vert \sum_{j=K_n}^{n-1} X_{j-l} \e_{j+1}\right \vert^q \right] \leq C \left(\frac{1}{n-K_n+1} \sum_{s=K_n}^{n-1} \sum_{t=K_n}^{n-1} \sigma(s-t) \sigma_\e(s-t)\right)^{q/2} 
\end{displaymath}
where $\sigma_\e(.)$ is the autocovariance function of the process $(\e_t)_{t \in \mathbb{Z}}$ i.e.
\begin{displaymath}
\sigma_\e(s-t)= \E(\e_t \e_s) =
\left\lbrace \begin{array}{l l}
 0 & \textrm{if}\,s\neq t\\
 1&\textrm{otherwise}
\end{array}\right. 
.
\end{displaymath}
We obtain:
\begin{eqnarray*}
\E\left[ (n-K_n+1)^{-q/2} \left \vert \sum_{j=K_n}^{n-1} X_{j-l} \e_{j+1}\right \vert^q \right] &\leq& C \left( \frac{1}{n-K_n+1}(n-K_n+1)\sigma(0)\right)^{q/2} \\
 &=& \OO(1).
\end{eqnarray*}
That concludes the proof. 
\end{proof}

\subsection{Proof of Lemma \ref{lmcalc}}
We recall that the constant $C$ may have different values in the different equations but is always independent of $n$ and then of $k$ since we want a convergence for all $1 \leq k \leq K_n$.\\
We decompose $f(k)- f_1(k)$ into 4 parts, which we estimate separately:
\begin{eqnarray}
 f(k)- f_1(k)&=& \textbf{X}_n^{*'} (k)\Sigma^{-1}(k)\frac{1}{n-K_n+1}\left( \sum_{j=K_n}^{n-\sqrt{n}-1}\textbf{X}_j(k)\e_{j+1}-\sum_{j=K_n}^{n-1}\textbf{X}_j(k)\e_{j+1}\right) \label{decomp1}\\
 &&+\textbf{X}_n^{*'} (k)\left( \Sigma^{-1}(k)-\widehat{\Sigma}^{-1}_n(k)\right) \frac{1}{n-K_n+1}\sum_{j=K_n}^{n-1}\textbf{X}_j(k)\e_{j+1}\label{decomp2}\\
 &&+\left( \textbf{X}_n^{*'}-\textbf{X}_n'(k)\right) \widehat{\Sigma}^{-1}_n(k)\frac{1}{n-K_n+1}\sum_{j=K_n}^{n-1}\textbf{X}_j(k)\e_{j+1}\label{decomp3}\\
 &&+\textbf{X}_n'(k)\widehat{\Sigma}^{-1}_n(k)\frac{1}{n-K_n+1}\sum_{j=K_n}^{n-1}\textbf{X}_j(k)\left(\e_{j+1} -\e_{j+1,k}\right) \label{decomp4}
\end{eqnarray}

\subparagraph*{Study of the term given in \eqref{decomp1}}

In this part we want to prove the mean-squared convergence to 0 of:
\begin{eqnarray}
&&\sqrt{\frac{n-K_n+1}{k}}\textbf{X}_n^{*'} (k)\Sigma^{-1}(k)\frac{1}{n-K_n+1}\left[ \sum_{j=K_n}^{n-\sqrt{n}-1}\textbf{X}_j(k)\e_{j+1}-\sum_{j=K_n}^{n-1}\textbf{X}_j(k)\e_{j+1}\right]\label{term1}\\ &= &
 \sqrt{\frac{1}{k}}\textbf{X}_n^{*'}\Sigma^{-1}(k)\frac{1}{\sqrt{n-K_n+1}}\sum_{j=n-\sqrt{n}-1}^{n-1}\textbf{X}_j(k)\e_{j+1}. \nonumber
\end{eqnarray}
Hölder's inequality applied twice with $1/p+1/q=1$ and $1/p'+1/q'=1$ gives:
\begin{eqnarray*}
 &&\E\left(\sqrt{\frac{1}{k}}\textbf{X}_n^{*'}(k)\Sigma^{-1}(k)\frac{1}{\sqrt{n-K_n+1}}\sum_{j=n-\sqrt{n}-1}^{n-1}\textbf{X}_j(k)\e_{j+1} \right)^2 \\ 
&\leq&\left( \E \left\Vert\frac{1}{\sqrt{k}}\textbf{X}_n^{*'}(k)\Sigma^{-1}(k)\right\Vert^{2q} \right)^{1/q} \left( \E \left \Vert\frac{1}{\sqrt{n-K_n+1}}\sum_{j=n-\sqrt{n}-1}^{n-1}\textbf{X}_j(k)\e_{j+1} \right \Vert^{2p} \right)^{1/p} \\
&\leq&\left( \E \left\Vert\frac{1}{\sqrt{k}}\textbf{X}_n^{*'}(k) \right \Vert^{2q'q}\right)^{1/(q'q)}\left( \E \left\Vert\Sigma^{-1}(k) \right \Vert^{2p'q}\right)^{1/(p'q)}\left( \E \left \Vert\frac{1}{\sqrt{n-K_n+1}}\sum_{j=n-\sqrt{n}-1}^{n-1}\textbf{X}_j(k)\e_{j+1} \right \Vert^{2p} \right)^{1/p} .
\end{eqnarray*}
Under assumption H.3 for all $p'$ and $q'$:
\begin{equation}
 \left( \E \left\Vert\Sigma^{-1}(k) \right \Vert^{2p'q}\right)^{1/(p'q)}=\left \Vert\Sigma^{-1}(k) \right \Vert^{2}=\OO(1) \label{maj1.1}
\end{equation}
since the spectral density of the process $(X_t)_{t \in \mathbb{Z}}$ admits a positive lower bound and then the largest eigenvalue of $\Sigma^{-1}(k) $ is bounded.
Furthermore by the convexity of the function $x \mapsto x^{2q'q}$ and the stationarity of the process $(\e_t)_{t \in \mathbb{Z}}$,
\begin{eqnarray*}
 \left( \E \left\Vert\frac{1}{\sqrt{k}}\textbf{X}_n^{*'}(k) \right \Vert^{2q'q}\right)^{1/(q'q)}&\leq&\left( \E   \left[\sum_{j=0}^{\sqrt{n}/2-K_n} b_j\e_{n-j} \right]^{2q'q} \right)^{1/(q'q)}
\end{eqnarray*}
and by Lemma 2 of \cite{wei}:
\begin{eqnarray}
 \left( \E \left\Vert\frac{1}{\sqrt{k}}\textbf{X}_n^{*'}(k) \right \Vert^{2q'q}\right)^{1/(q'q)}
 &\leq& C' \left(\sum_{j=0}^{\sqrt{n}/2-K_n}b_j^2 \right)  \nonumber \\
 & \leq &C\label{x*}
\end{eqnarray}
 because the sequence $(b_j^2)_{j \in \mathbb{N}}$ is summable.\\
Finally by Lemma \ref{lem2}:
\begin{eqnarray}
 &&\left( \E \left \Vert\frac{1}{\sqrt{n-K_n+1}}\sum_{j=n-\sqrt{n}-1}^{n-1}\textbf{X}_j(k)\e_{j+1} \right \Vert^{2p} \right)^{1/p}\nonumber \\
 & \leq&\frac{(n+1)^{1/4}}{\sqrt{n-K_n}}\left( \E \left \Vert\frac{1}{(n+1)^{1/4}}\sum_{j=n-\sqrt{n}-1}^{n-1}\textbf{X}_j(k)\e_{j+1} \right \Vert^{2p} \right)^{1/p} \nonumber\\
&\leq & C\left(\frac{k}{n^{1/4}} \right) \label{maj1.3}
\end{eqnarray}
By inequalities \eqref{maj1.1}, \eqref{x*} and \eqref{maj1.3}:
\begin{displaymath}
\E\left(\sqrt{\frac{1}{k}}\textbf{X}_n^{*'}(k)\Sigma^{-1}(k)\frac{1}{\sqrt{n-K_n+1}}\sum_{j=n-\sqrt{n}-1}^{n-1}\textbf{X}_j(k)\e_{j+1} \right)^2  \leq C \left(\frac{k}{n^{1/4}}\right) \leq C \left(\frac{K_n}{n^{1/4}} \right) 
\end{displaymath}
which converges to 0 as $n$ tends to infinity under assumption \eqref{condth2}.

\subparagraph*{Study of the term given in \eqref{decomp2}}

Prove that: \begin{equation}
 \lim_{n \rightarrow + \infty}\E\left( \sqrt{\frac{n-K_n+1}{k}}\textbf{X}_n^{*'} (k)\left[\Sigma^{-1}(k)- \widehat{\Sigma}^{-1}_n(k)\right] \sum_{j=K_n}^{n-1}\textbf{X}_j(k)\e_{j+1}\right)^2 =0. \label{term2}
\end{equation}
Applying twice Hölder's inequality, we have:
\begin{eqnarray*}
 &&\E\left( \sqrt{\frac{n-K_n+1}{k}}\textbf{X}_n^{*'} (k)\left[\Sigma^{-1}(k)- \widehat{\Sigma}^{-1}_n(k)\right] \sum_{j=K_n}^{n-1}\textbf{X}_j(k)\e_{j+1}\right)^2 \\
& \leq &  \left( \E \left\Vert\frac{1}{\sqrt{k}}\textbf{X}_n^{*'}(k) \right \Vert^{2q'q}\right)^{1/(q'q)}\left( \E \left\Vert\Sigma^{-1}(k)-\widehat{\Sigma}^{-1}_n(k) \right \Vert^{2p'q}\right)^{1/(p'q)}\\&&\left( \E \left \Vert\frac{1}{\sqrt{n-K_n+1}}\sum_{j=K_n}^{n-1}\textbf{X}_j(k)\e_{j+1} \right \Vert^{2p} \right)^{1/p} .
\end{eqnarray*}
Applying Lemma \ref{lem2} we obtain that:
\begin{equation}
 \left( \E \left \Vert\frac{1}{\sqrt{n-K_n+1}}\sum_{j=K_n}^{n-1}\textbf{X}_j(k)\e_{j+1} \right \Vert^{2p} \right)^{1/p}\leq C k\leq C K_n  \label{maj2.1}
\end{equation}
since $k\leq K_n$.
Now we derive the mean-squared convergence to 0 when $d \in ]0,1/2[$.
\par For $d \in]0,1/4[$, we apply Theorem \ref{th1} and we get:
\begin{equation}\label{maj2.2}
 \left( \E \left\Vert\Sigma^{-1}(k)-\widehat{\Sigma}^{-1}_n(k) \right \Vert^{2p'q}\right)^{1/(p'q)}\leq C\left( \frac{K_n^2}{n-K_n+1}\right) .
\end{equation}
Then it follows from \eqref{x*}, \eqref{maj2.1} and \eqref{maj2.2} that:
\begin{displaymath}
 \E\left( \sqrt{\frac{n-K_n+1}{k}}\textbf{X}_n^{*'} (k)\left[\Sigma^{-1}(k)- \widehat{\Sigma}^{-1}_n(k)\right] \sum_{j=K_n}^{n-1}\textbf{X}_j(k)\e_{j+1}\right)^2 \leq C\left( \frac{K_n^3}{n-K_n+1}\right) 
\end{displaymath}
which converges to 0 if condition \eqref{condth2} holds.
\par If $d \in]1/4,1/2[$, we obtain by Theorem \ref{th1} that
\begin{equation}
 \left( \E \left\Vert\Sigma^{-1}(k)-\widehat{\Sigma}^{-1}_n(k) \right \Vert^{2p'q}\right)^{1/(p'q)}\leq C \left( \frac{K_n^2 L^2(n-K_n+1)}{(n-K_n+1)^{2-4d}}\right) .\label{maj2.3}
\end{equation}
The inequalities \eqref{x*},\eqref{maj2.1} and \eqref{maj2.3} allow us to conclude that:
\begin{displaymath}
 \E\left( \sqrt{\frac{n-K_n+1}{k}}\textbf{X}_n^{*'} (k)\left[\Sigma^{-1}(k)- \widehat{\Sigma}^{-1}_n(k)\right] \sum_{j=K_n}^{n-1}\textbf{X}_j(k)\e_{j+1}\right)^2 \leq C \left( \frac{K_n^3 L^2(n-K_n+1)}{(n-K_n+1)^{2-4d}}\right) 
\end{displaymath}
which converges to 0 under assumption \eqref{condth2}.
\par For $d=1/4$, applying Theorem \eqref{th1} we obtain that
\begin{equation}\label{maj2.4}
 \left( \E \left\Vert\Sigma^{-1}(k)-\widehat{\Sigma}^{-1}_n(k) \right \Vert^{2p'q}\right)^{1/(p'q)}\leq C \left( \frac{K_n^2 L^2(n-K_n+1)\log(n-K_n+1)}{(n-K_n+1)}\right).
\end{equation}
The inequalities \eqref{x*}, \eqref{maj2.1} and and \eqref{maj2.4} allow us to conclude that:
\begin{displaymath}
 \E\left( \sqrt{\frac{n-K_n+1}{k}}\textbf{X}_n^{*'} (k)\left[\Sigma^{-1}(k)- \widehat{\Sigma}^{-1}_n(k)\right] \sum_{j=K_n}^{n-1}\textbf{X}_j(k)\e_{j+1}\right)^2 \leq C \left( \frac{K_n^3 L^2(n-K_n+1)\log(n-K_n+1)}{(n-K_n+1)}\right) 
\end{displaymath}
\subparagraph*{Study of the term given in \eqref{decomp3}}

 Prove that:
\begin{equation}\label{term3}
 \lim_{n \rightarrow + \infty}\E\left( \sqrt{\frac{n-K_n+1}{k}}\left[ \textbf{X}_n^{*'} (k)-\textbf{X}_n'(k)\right] \widehat{\Sigma}^{-1}_n(k) \sum_{j=K_n}^{n-1}\textbf{X}_j(k)\e_{j+1}\right)^2 =0.
\end{equation}
Using Holder's inequality twice, we have: 
\begin{eqnarray}
 &&\E\left( \sqrt{\frac{n-K_n+1}{k}}\left[ \textbf{X}_n^{*'} (k)-\textbf{X}_n'(k)\right]  \widehat{\Sigma}^{-1}_n(k) \sum_{j=K_n}^{n-1}\textbf{X}_j(k)\e_{j+1}\right)^2  \nonumber \\
&\leq& \left( \E \left\Vert\frac{1}{\sqrt{k}}\left[ \textbf{X}_n^{*'}(k) -\textbf{X}_n'(k)\right]  \right \Vert^{2q'q}\right)^{1/(q'q)}\left( \E \left\Vert \widehat{\Sigma}^{-1}_n(k) \right \Vert^{2p'q}\right)^{1/(p'q)}\nonumber\\
&&\left( \E \left \Vert\frac{1}{\sqrt{n-K_n+1}}\sum_{j=K_n}^{n-1}\textbf{X}_j(k)\e_{j+1} \right \Vert^{2p} \right)^{1/p} .\nonumber
\end{eqnarray}
In view of the convexity of the function $x \mapsto x^{qq'}$ and of the stationarity of the process $(\e_t)_{t \in \mathbb{Z}}$, we have:
\begin{displaymath}
 \max_{1 \leq k \leq K_n}\left( \E \left\Vert\frac{1}{\sqrt{k}}\left[\textbf{X}_n^{*'}(k) -\textbf{X}_n'(k) \right] \right \Vert^{2q'q}\right)^{1/(q'q)} \leq \left( \E \left(\sum_{j=\sqrt{n}/2-K_n+1}^{+\infty}b_j \e_{n-j-l}  
 \right )^{2q'q}\right)^{1/(q'q)}. 
\end{displaymath}
And by Lemma 2 of \cite{wei}, we obtain
\begin{eqnarray}
 \max_{1 \leq k \leq K_n}\left( \E \left\Vert\frac{1}{\sqrt{k}}\left[\textbf{X}_n^{*'}(k) -\textbf{X}_n'(k) \right] \right \Vert^{2q'q}\right)^{1/(q'q)}
 & \leq &C\left(\sum_{j=\sqrt{n}/2-K_n+1}^{+\infty} b_j^2 \right)\nonumber \\
& \leq &C n^{\frac{2d-1}{2}}.\label{maj3.1}
\end{eqnarray}
Then by Theorem \ref{th1}:
\begin{equation}
 \left( \E \left\Vert\widehat{\Sigma}^{-1}_n(k) \right \Vert^{2p'q}\right)^{1/(p'q)}\leq C. \label{maj3.2}
\end{equation}
By inequalities \eqref{maj2.1}, \eqref{maj3.1} and \eqref{maj3.2}, we then obtain:
\begin{displaymath}
 \E\left( \sqrt{\frac{n-K_n+1}{k}}\left[ \textbf{X}_n^{*'} (k)-\textbf{X}_n'(k)\right]  \widehat{\Sigma}^{-1}_n(k) \sum_{j=K_n}^{n-1}\textbf{X}_j(k)\e_{j+1}\right)^2\leq C\left( K_n n^{\frac{2d-1}{2}}\right) 
\end{displaymath}
which converges to 0 as $n$ tends to infinity if condition \eqref{condth2} holds.

\subparagraph*{Study of the term given in \eqref{decomp4}}

We want to prove that:
\begin{equation}
 \lim_{n \rightarrow + \infty} \E\left( \sqrt{\frac{1}{\E\left( \mathcal{S}_n(k)\right)^2}}\textbf{X}_n' (k)\widehat{\Sigma}^{-1}_n(k)\frac{1}{n-K_n+1}\sum_{j=K_n}^{n-1}\textbf{X}_j(k)\left[ \e_{j+1,k}-\e_{j+1}\right] \right)^2 =0.\label{term4}
\end{equation}
Using Hölder's inequality twice, we have:
\begin{eqnarray*}
&& \E\left( \sqrt{\frac{1}{\E\left( \mathcal{S}_n(k)\right)^2}}\textbf{X}_n' (k)\widehat{\Sigma}^{-1}_n(k)\frac{1}{n-K_n+1}\sum_{j=K_n}^{n-1}\textbf{X}_j(k)\left[ \e_{j+1,k}-\e_{j+1}\right] \right)^2\\
&\leq& \frac{1}{(n-K_n+1)\E\left( \mathcal{S}_n(k)\right)^2}\left( \E \left\Vert\frac{1}{\sqrt{n-K_n+1}}\sum_{j=K_n}^{n-1}\textbf{X}_j(k)\left[ \e_{j+1,k}-\e_{j+1}\right]  \right \Vert^{2q'q}\right)^{1/(q'q)}\\  
&&\left( \E \left\Vert\widehat{\Sigma}^{-1}_n(k) \right \Vert^{2p'q}\right)^{1/(p'q)}\left( \E \left \Vert\textbf{X}_n'\right \Vert^{2p} \right)^{1/p} .
\end{eqnarray*}
Applying Lemma \ref{lem1}, we obtain for every $\delta>0$:
\begin{equation}
\left( \E \left\Vert\frac{1}{\sqrt{n-K_n+1}}\sum_{j=K_n}^{n-1}\textbf{X}_j(k)\left[ \e_{j+1,k}-\e_{j+1}\right]  \right \Vert^{2q'q}\right)^{1/(q'q)} \leq C\left(  k (n-K_n+1)^{2d+\delta}\E\left( \mathcal{S}_n(k)\right)^2\right) .\label{maj4.1}
\end{equation}
Finally we choose $p=2$ and we have:
\begin{eqnarray*}
\left( \E \left\Vert\textbf{X}_n'(k)\right \Vert^{4} \right)^{1/2}&=& \sqrt{\E\left(\sum_{j=1}^kX_j^2 \right)^2 } \\
&=& \sqrt{\left(\sum_{j=1}^k\sigma(j)\right)^2+2\sum_{j=1}^k\sum_{l=1}^k\sigma(j-l)^2}
\end{eqnarray*}
since the process $(X_t)_{t \in \mathbb{Z}}$ is Gaussian. Using the assumption H.2 on the covariances, we verify that for all $\delta>0$:
\begin{equation}
 \left( \E \left\Vert\textbf{X}_n'(k)\right \Vert^{4} \right)^{1/2} \leq C\sqrt{ k^{4d+\delta}} 
\leq C \sqrt{k} \label{momx}
\end{equation}
if $d \in \left]  0,1/4 \right[ $.
With these three inequalities \eqref{maj3.2}, \eqref{maj4.1}, \eqref{momx} and $1\leq k \leq K_n$, we conclude that:
\begin{eqnarray*}
 \forall \delta>0, &&\: \E\left( \sqrt{\frac{1}{\E\left( \mathcal{S}_n(k)\right)^2}}\textbf{X}_n' (k)\widehat{\Sigma}^{-1}_n(k)\frac{1}{n-K_n+1}\sum_{j=K_n}^{n-1}\textbf{X}_j(k)\left( \e_{j+1,k}-\e_{j+1}\right) \right)^2\\
& \leq & C\frac{1}{(n-K_n+1)\E\left( \mathcal{S}_n(k)\right)^2}\left( \sqrt{k} k (n-K_n+1)^{2d+ \delta}\E\left( \mathcal{S}_n(k)\right)^2\right) \\
& \leq & C\frac{K_n^{3/2}}{(n-K_n+1)^{1-2d-\delta}}.
\end{eqnarray*}
which converges to 0 under condition \eqref{condth2}.\\
On the other hand if $d \in [1/4,1/2[$, inequality \eqref{momx} becomes:
\begin{equation}
 \forall \delta >0,\: \left( \E \left\Vert\textbf{X}_n'(k)\right \Vert^{4} \right)^{1/2}\leq C k^{4d+ \delta} \label{maj4.2}
\end{equation}
Using  inequalities \eqref{maj3.2}, \eqref{maj4.1} and \eqref{maj4.2}, we have for all $\delta>0$
\begin{displaymath}
\E\left( \sqrt{\frac{1}{\E\left( \mathcal{S}_n(k)\right)^2}}\textbf{X}_n' (k)\widehat{\Sigma}^{-1}_n(k)\frac{1}{n-K_n+1}\sum_{j=K_n}^{n-1}\textbf{X}_j(k)\left[ \e_{j+1,k}-\e_{j+1}\right] \right)^2\leq  C\left( \frac{K_n^{1+2d+\delta}}{(n-K_n+1)^{1-2d-\delta}}\right) 
\end{displaymath}
which converges to 0 under condition \eqref{condth2}.

\par We have proved that for all $d \in \left]0, \frac{1}{2} \right[ $ \begin{displaymath}
\lim_{n \rightarrow + \infty} \max_{1 \leq k \leq K_n}\E\left(\sqrt{\frac{1}{L_n(k)}}(f(k)-f_1(k)) \right)^2 =0.
\end{displaymath}

\bibliographystyle{apalike} 
\bibliography{biblioprev}

\end{document}